\documentclass[11pt]{article}
\input{amssym.def}
\input{amssym}
\usepackage{eufrak}

\newtheorem{theorem}{Theorem}

\newtheorem{conjecture}[theorem]{Conjecture}
\newtheorem{corollary}[theorem]{Corollary}

\newtheorem{definition}[theorem]{Definition}

\newtheorem{proposition}[theorem]{Proposition}
\newtheorem{remark}[theorem]{Remark}

\newlength{\coefkappa}
\begin{document}

\def\h{{\hbar}}
\def\Cas{{\bf Cas}}
\def\Op{{\rm Op}}
\def\Caskm{{\bf Cas}_{(k,m)}}
\def\ch{{\cal CH}}
\def\chc{{\cal CH^{\chi}}}
\def\lhq{{\cal L}_{\h, q}}
\def\slhq{{\cal SL}_{\h, q}}
\def\lc{{\cal L}_{\h, q}^{\chi}}
\def\gggg{{\bf g}}
\def\gh{\gggg_{\h}}
\def\sll{s_{\la}}
\def\glh{gl(2)_{\h}}
\def\gq{\gggg_{q}}
\def\gqq{\gggg_{q}^{\ot 2}}
\def\Ugl{U(gl(V))}
\def\Usl{U(sl(V))}
\def\Ush{U(sl(2)_{\h})}
\def\Zgh{Z[U(gl(2)_{\h})]}
\def\Zsh{Z(U(sl(2)_{\h}))}
\def\uq{U_q(sl(2))}
\def\Uq{U_q(sl(n))}
\def\ugh{U(\gggg_{\h})}
\def\uqs{U_q(sl(n))}
\def\ug{U(\gggg)}
\def\C{{\Bbb C}}
\def\B{{\cal C}}
\def\bl{{B_{\la}}}
\def\cl{{C_{\la}}}
\def\cp{{\Bbb CP}}
\def\R{{\Bbb R}}
\def\K{{\Bbb K}}
\def\Z{{\Bbb Z}}
\def\ii{{\sqrt{-1}}}
\def\O{{\cal O}}
\def\A{{\cal A}}
\def\Ah{{\cal A}_\hbar}
\def\kO{{\K(\cal O)}}
\def\De{\Delta}
\def\de{\delta}
\def\la{\lambda}
\def\La{\Lambda}
\def\co{{\cal O}}

\def\Pp{P_+}
\def\Pm{P_-}
\def\Ppm{P_{\pm}}

\def\lp{\Lambda_+(V)}
\def\lm{\Lambda_-(V)}
\def\lpm{\Lambda_{\pm}(V)}
\def\lpmk{\Lambda_{\pm}^k(V)}
\def\lmk{\Lambda_{-}^k(V)}
\def\lmpp{\Lambda_{-}^{p-1}(V)}
\def\lmp{\Lambda_{-}^p(V)}

\def\l{{\bf l}}
\def\opi{\bar \pi}
\def\lpml{\Lambda_{\pm}^l(V)}
\def\lmp{\Lambda_{-}^p(V)}
\def\qm{q^{-1}}
\def\ve{\varepsilon}
\def\ep{\epsilon}
\def\al{{\alpha}}
\def\ah{{ A}_{\h}}
\def\aq{{ A}_{q}}
\def\ahq{{ A}_{\h,q}}
\def\ot{\otimes}
\def\om{\omega}
\def\Om{\Omega}
\def\Ml{M^{\lambda}}
\def\End{{\rm End\, }}
\def\Mat{{\rm Mat }}
\def\Mor{{\rm Mor\, }}
\def\tr{{\rm tr}}
\def\Tr{{\rm Tr}}
\def\Ind{{\rm Ind\, }}
\def\prePic{{\rm prePic }}
\def\Mod{{\rm Mod\, }}
\def\Uqmod{U_q(sl(n))-{\rm Mod}}
\def\uqmod{U_q(sl(2))-{\rm Mod}}
\def\slm{sl(n)-{\rm Mod}}
\def\Hom{{\rm End\, }}
\def\Sym{{\rm Sym\, }}
\def\Ker{{\rm Ker\, }}
\def\Im{{\rm Im\, }}
\def\Ob{{\rm Ob\, }}

\def\id{{\rm id\, }}
\def\dim{{\rm dim}}
\def\ddim{{\rm dim}}
\def\span{{\rm span\, }}
\def\ttr{{\bf tr\, }}
\def\ad{{\rm ad\, }}
\def\cotr{{\rm cotr\, }}
\def\rk{{\rm rk\, }}
\def\Vect{{\rm Vect\,}}
\def\Fun{{\rm Fun\,}}
\def\vv{V^{\ot 2}}
\def\vvv{V^{\ot 3}}
\def\vl{V_{\la}}
\def\vm{V_{\mu}}
\def\lla{\bf{\la}}

\def\bea{\begin{eqnarray}}
\def\eea{\end{eqnarray}}
\def\be{\begin{equation}}
\def\ee{\end{equation}}
\def\nn{\nonumber}

\makeatletter
\renewcommand{\theequation}{{\thesection}.{\arabic{equation}}}
\@addtoreset{equation}{section} \makeatother

\title{q-Index on braided  non-commutative spheres}
\author{D.~Gurevich, R.~Leclercq\\
{\small\it ISTV, Universit\'e de Valenciennes, 59304 Valenciennes,
France}\\
\rule{0pt}{7mm}P.~Saponov,\\
{\small\it Theory Department of Institute for High Energy
Physics, 142281 Protvino, Russia}}

\maketitle

\begin{abstract}
To some Hecke symmetries (i.e. Yang-Baxter braidings of
Hecke type) we assign algebras called braided
non-commutative spheres. For any such algebra, we introduce
and compute a q-analog of the Chern-Connes index. Unlike
the standard Chern-Connes index, ours is based on the
so-called categorical trace specific for a braided category
in which the algebra in question is represented.
\end{abstract}

{\bf AMS Mathematics Subject Classification, 1991:} 17B37, 81R50


\section{Introduction}
\label{sec:int}

The Chern-Connes index (pairing) was introduced  and
intensively used for needs of non-com\-mu\-ta\-tive
geometry by A. Connes. Given an associative algebra $A$,
the Chern-Connes index (in its simplest form, cf. \cite{L})
can be defined via a pairing \be \Ind:K_0(A)\times
K^0(A)\to \Z,\label{ind} \ee where $K^0(A)$ is the
Grothendieck group of the monoid of its finite dimensional
representations and $K_0(A)$ is the Grothendieck group of
classes of one-sided  projective
$A$-modules\footnote{Throughout the paper all projective
modules are supposed to be finitely generated.}. (Observe
that according to the Serre-Swan approach such modules are
considered as appropriate analogs of vector bundles on a
variety.) We will only deal with the finite dimensional
representations of algebras in question taking the category
of finite dimensional $U(sl(n))$-modules as a pattern. In
this sense our setting is purely algebraic. This is the
main difference between our approach and that based on the
Connes spectral triples in which a considerable amount of
functional analysis is involved (cf. \cite{C}, where the
function algebra on the quantum group $\uq$ is studied from
this viewpoint).

Any projective module can be identified  with an idempotent
$e\in\Mat(A)$ where as usual $\Mat(A) =\oplus_n \Mat_n(A)$
and $\Mat_n(A)$ stands for the algebra of $n\times n$
matrices with entries from $A$.

Let us fix a representation $\pi_U:A\to \End(U)$ and an
idempotent $e\in \Mat(A)$.  Then the pairing (\ref{ind}) is defined by
\be
\Ind(e,\pi_U)= \tr(\pi_U(\tr\, e))=\tr(\pi_U(e)), \label{Ind}
\ee
where $\pi_U$ is   naturally extended to $\Mat(A)$.

On defining  $K_0(A)$ in the standard way (cf. \cite{Ro})
it is not difficult to show that $\Ind(e,\pi_U)$ does not
depend on a representative of a class from $K_0(A)$ and
thus pairing (\ref{Ind}) reduces to (\ref{ind}).

In this paper we introduce a ``braided version" of the Chern-Connes
index.
This
version is based on the so-called ``categorical trace" (see section
\ref{sec:cat}) and motivated
by the ``braided nature" of the algebras considered. These algebras are
quotients of some braided
analogs of enveloping algebras $U(gl(n))$ and $U(sl(n))$ and are thought
of as braided
non-commutative counterparts of orbits in $sl(n)^*$.
As a result we compute our index on a class of orbits of this type
called
``braided spheres", in  particular, for the ``quantum sphere" related to
the quantum group $\uq$.

First of all, let us briefly
describe the braided categories we are working with. The detailed
consideration
is presented in section \ref{sec:cat}. Any such
category is generated by a finite dimensional vector space $V$
equipped with a map to be called {\em a braiding (morphism)}
\be
R:\vv\to \vv \label{brdn}
\ee
which satisfies the quantum Yang-Baxter equation
\be
R_{12}R_{23}R_{12}=R_{23}R_{12}R_{23},\quad R_{12}=R\ot\id, \quad
R_{23}=\id\ot R.\label{QYBE}
\ee
Besides, we will suppose $R$ to be of the Hecke type. This means
that the braiding $R$ satisfies the following {\em Hecke
condition}
\be
(q\,\id-R)(q^{-1}\id +R)=0,\quad q\in \K.
\label{He-cond}
\ee
Hereafter $\K$ stands for the basic field (usually $\C$ but sometimes
$\R$ is allowed) and the
parameter $q\in\K$ is assumed to
be generic (but $q=1$ is permitted). The braidings of the Hecke type
will be also called
{\em Hecke symmetries}.

Let $\B=\B(V)$ be the category generated by the space $V$. The
sets of its objects and categorical morphisms will be denoted
respectively by $\Ob(\B)$ and $\Mor(\B)$. The category $\Uqmod$ of all
finite
dimensional modules over the quantum group $\Uq$ serves as an example
of $\B(V)$.
In this case the space $V$ is the fundamental (vector) module, the
braiding $R$ is the
Drinfeld-Jimbo $R$-matrix and the categorical morphisms are
linear maps commuting with the action of the quantum group $\Uq$.

Under some additional conditions on $R$ (see section \ref{sec:cat})
the braided categories in
question are {\it rigid} (for the terminology the reader is referred to
\cite{CP}). This means that for any $U\in \Ob(\B)$ there exists
$U^*_r\in\Ob(\B)$
(resp. $U^*_l\in\Ob(\B)$) for which one can define a non-degenerate
pairing
$$
U\ot U^*_r\to \K \quad ({\rm resp.}\;\,U^*_l\ot U\to \K)
$$
and this map is a categorical morphism. The space $U^*_r$
(resp. $U^*_l$) is called the {\em right (left) dual} space
to $U$. Therefore, for any $U\in\Ob(\B)$ the space
of its right (resp. left) ``internal endomorphisms" 
\be
\End_r(U)\stackrel{\mbox{\tiny def}}{=} U_r^*\ot
U,\quad({\rm resp.}\quad \End_l(U)\stackrel{\mbox{\tiny
def}}{=} U\ot U_l^*)\label{endo} 
\ee 
is also contained in $\Ob(\B)$.

Then, in $\B$ we define an important categorical morphism
$$
\tr_R:\quad \End_\varepsilon (U)\to\K, \quad {\varepsilon = l,r}
$$
called {\em the categorical trace}. The super-trace
is an example of such a categorical trace. Namely, in super-algebra and
super-geometry
this trace replaces the classical one. For a similar reason, dealing
with a braided
category, we make use of the categorical trace specific for this
category.

Now, let us pass to algebras in question. Assume for a moment that
$q=1$. This means that
our braiding $R$ becomes involutive: $R^2=\id$. For such a braiding
there exists a natural way
to define {\em a generalized Lie bracket}:
\be \
[\,\,,\,\,]:\quad \End_\varepsilon(V)^{\ot 2}\to \End_\varepsilon(V),
 \quad {\varepsilon = l,r}\label{Liebr}
\ee
(cf. \cite{G} for detail). Being equipped with such a bracket, the
space $\End_l(V)$
(for definiteness we set $\varepsilon=l$) becomes {\em a generalized
Lie algebra}. It
will be denoted  $gl_R(V)$. For instance, a super-Lie algebra is a
particular case of generalized one.

Moreover, for the aforementioned categorical trace the
subspace $sl_R(V)$ of all traceless elements is closed with
respect to this bracket. Thus, the space $sl_R(V)$ is also
{\em a generalized Lie algebra}. Then their enveloping
algebras $U(gl_R(V))$ and $U(sl_R(V))$ can be defined by
systems of quadratic-linear equations. Furthermore, they
become {\em braided Hopf algebras}, being equipped with an
appropriate coproduct, antipode and counit. On generators
$X\in gl_R(V)$ (or $X\in sl_R(V)$) this coproduct has the
classical form: \be \Delta(X)=X\ot 1 + 1\ot X.\label{copr}
\ee

By means of the coproduct (which gives rise to a braided version
of the Leibniz rule) we can construct an embedding
$$
\End_l(V)\to\End_l(V^{\ot m}),\,\,\forall m.
$$
Restricting these maps to the subspaces associated with the
Young diagrams (the corresponding Young projectors can be
constructed for any involutive braiding) we get a family of
irreducible representations of the generalized Lie algebra
$sl_R(V)$. Then considering all their direct sums we get a
category of finite dimensional representations of the
algebra in question  similar to that of $sl(n)$, $n=\dim
V$. The main difference between the later category and a
braided one consists in traces. When considering the
algebra $U(gl_R(V))$ (or $U(sl_R(V))$) with the
aforementioned category of finite dimensional
representations it is natural to use the corresponding
categorical trace  in order to define all numerical
characteristics (dimensions, indexes, etc).

It is this scheme that is realized in the present paper.
However, here we deal with a more interesting (and more
difficult) case of algebras and categories associated to
some non-involutive Hecke symmetries. In this case it is
not evident which algebras should be taken as $U(gl_R(V))$
and $U(sl_R(V))$. The matter is that a direct
generalization of bracket (\ref{Liebr}) to a non-involutive
$R$ leads to ``enveloping algebras'' which are not flat
deformations of the classical ones even if $R$ is a
deformation of the usual flip. Otherwise stated, the
dimensions of the homogeneous components of the
corresponding graded algebra differ from their classical
analogs.

Nevertheless, there exist algebras (denoted below as $\lhq$ and $\slhq$)
possessing good
deformational properties and playing the role of braided analogs of
enveloping algebras $U(gl(n))$ and
$U(sl(n))$ respectively (though apparently they do not look like the
usual enveloping algebras). They can
be described in terms of the {\it modified reflection equation} (mRE).
This equation can be defined for
any braiding;  for involutive $R$ it leads to the aforementioned
enveloping algebras $U(gl_R(V))$ and
$U(sl_R(V))$. The algebras $\lhq$ and $\slhq$ are described in the next
section.

In the case related to the quantum group $\Uq$ (in the $\Uq$ case for
short)
these algebras are one-sided
$\Uq$--modules unlike the $\Uq$ themself which is a two-sided
$\Uq$--module.
Using the results of \cite{LS}
it is possible to show that any $\Uq$--module becomes an $\slhq$--one.
For the Hecke symmetry coming from
$\Uq$ the corresponding category $\B(V)$ is the representation category
of this quantum group and all its
objects can be equipped with an action of the algebra  $\slhq$. The
family 
of all representations of the
algebra $\slhq$ is, however, larger than that of $\Uq$. Nevertheless, we
only need finite dimensional
$\slhq$--representations which are $\Uq$--modules and are equivariant
(covariant) with respect to
the action of the quantum group. (The problem of describing other
representations of the algebras in question
will be considered elsewhere. Here we only consider the representation
theory in the category $\B(V)$.)

But in general case we do not have such a useful tool as the quantum
group.
So, we modify the notion of equivariant
representation in order to adapt it to a new setting.

Futhermore, restricting ourselves to the case $\rk(R)=2$ (see section
\ref{sec:cat}) we equip the category $\B(V)$
with an equivariant action of the algebra $\slhq$. Then we introduce
``a braided non-com\-mu\-ta\-tive sphere" which
is a quotient of this algebra (we get it by fixing a value of the
quadratic
braided Casimir element). In the $\uq$ case the
braided sphere is also called ``the quantum non-com\-mu\-ta\-tive
sphere''.

Our quantum sphere is close to the known Podles sphere. However, while
the
Podles sphere is
$U_q(su(2))$--homogeneous space and is introduced via some reduction
from
its dual, our quantum sphere
is defined via the mRE algebra as its approriate quotient. Consequently,
the representation theory of
the Podles quantum sphere constructed in \cite{P} differs drastically
from
that of $\slhq$.

The Chern-Connes index on the Podles sphere was computed in the work
\cite{H} with the use of the representation
theory from \cite{P}, the trace defined in  \cite{MNW}  and idempotents
introduced in \cite{HM}. These idempotenets
are labeled by $n\in\Z$ and are presented via generators of the function
algebra on the quantum group $U_q(su(2))$.

In contrast, our method of constructing projective modules over the
braided
(in particular, quantum) non-commutative
spheres makes use of a braided version of the Cayley-Hamilton identity
for
some matrices from $\Mat(\lhq)$. This
allows us to construct a larger set of projective modules. However, this
set
of modules turns out to be ``too large''
and this leads us to the problem of defining a reasonable equivalence
between modules in order to get the group $K_0$
of the classical size. This problem is discussed in the last section.

In the $\uq$ case all ingredients of our construction (representation
theory, projective modules, traces, indices)
have the classical limits at $q\to 1$. The corresponding algebra was
earlier
considered in \cite{GS2} under the name
of non-commutative sphere.

The paper is organized as follows. In the next section we
define the categories and algebras we are dealing with. In
section 3 we introduce and compute the braided version of
the Chern-Connes index. In section 4 we consider the
quantum sphere as an example of our general construction.
In section 5 some problems are discussed which arise in
connection with our  approach.

{\bf Acknowledgement} The authors would like to thank the referee for
valuble remarks. Two of us (D.G. and P.S.)
are grateful to Max-Planck-Institut f\"ur Mathematik (Bonn) where the
final version of the paper was written for
warm hospitality and stimulating atmosphere.

\section{Categories $\B(V)$ and related algebras}
\label{sec:cat}

We begin this section with a short description of the category
$\B=\B(V)$ generated by a finite dimensional vector space V
equipped with a Hecke symmetry $R$. This category forms a base
of all our considerations, for its detailed description see \cite{GLS}.

Given a Hecke symmetry $R$, one can connect with it a ``symmetric"
(resp. ``skew-symmetric") algebra
$\lp$ (resp. $\lm$) of the space $V$ defined as the quotient
$$
\lp=T(V)/\{\Im(q\,\id-R)\} \quad({\rm resp.}\;\;
\lm=T(V)/\{\Im(q^{-1}\,\id+R)\}).
$$
Here $T(V)$ stands for the free tensor algebra. Let $\lpmk$ be the
homogeneous component of $\lpm$
of degree $k$. If there exists an integer $p$ such that $\lmk$ is
trivial for $k>p$ and $\ddim(\lmp)=1$,
then $R$ is called {\em  an even symmetry} and $p$ is called {\em the
rank}
of $R$: $p=\rk(R)$. Hereafter
the symbol ``dim'' stands for the classical dimensions. In what follows
all
Hecke symmetries are
assumed to be even.

Using the Yang-Baxter equation (\ref{QYBE}) we can extend braiding
(\ref{brdn}) onto any tensor
powers of $V$
\be
R:\quad V^{\ot m}\ot V^{\ot n}\to V^{\ot n}\ot
V^{\ot m}\label{braid}
\ee
(as usual, we put $V^{\ot 0}=\K$ and $(1\ot x)\triangleleft
R = x\ot 1,\,\,(x\ot 1)\triangleleft R
=1\ot x, \,\,\forall x\in V^{\ot m}$).

For an arbitrary fixed integer $m\geq 2$ we consider partitions
$\lambda \vdash m$
$$
\la=(\la_1,\la_2,...,\la_k),\quad \la_1\ge\la_2\ge...\ge\la_k > 0,
\quad \la_1+...+\la_k=m
$$
($k$ is called {\em the height} of $\la$). There exists a natural
way to assign a space $\vl$ (equipped with a set of embeddings
$\vl\hookrightarrow V^{\ot m}$) to any partition $\la$ (cf.
\cite{GLS}). By definition the spaces $\vl$ are simple objects of
the category $\B$ (we will motivate this definition below). All
other objects are the direct sums of  the simple ones.

Let us describe {\it the categorical morphisms} in $\B$.
The term {\em categorical} emphasizes the difference among
morphisms from ${\rm Mor}(\B)$ and ``internal endomorphisms''
which are elements of $\End_\varepsilon(U)\in {\rm
Ob}(\C)$, $\epsilon=r,l$.

We distinguish categorical morphisms of two kinds.
Categorical morphisms of {\em the first kind} are the linear maps
$V^{\ot m}\to V^{\ot m}$ $m\ge 0$ coming from the Hecke
algebra, as well as their restrictions to any object embedded into
$V^{\ot m}$. Recall, that the Hecke algebra $H_m$ can be represented
in $V^{\ot m}$ by means of the Hecke symmetry $R$.

The categorical morphisms {\em of the second kind} arise from a
procedure of cancelling columns of height $p$ in the Young
diagram corresponding to a given partition $\lambda$. The procedure
is as follows. Denote by $v$ a generator of the one
dimensional space $\lmp$ (we call $v$ {\em the associated
determinant}). If one fixes a basis $\{x_i\}$ in the space $V$,
then $v$ can be represented in the form
$$
v=v^{i_1\dots i_p}x_{i_1}\otimes \dots \otimes x_{i_p}.
$$
Hereafter the summation over repeated indices is always
understood. The tensor $v^{i_1\dots i_p}$ is one of the two
structure tensors that define the matrix of the highest order
antisymmetrizer (projector) $A^{(p)}: \,V^{\otimes p}\rightarrow
\Lambda^p_-(V)$ in the basis $x_{i_1}\otimes \dots \otimes
x_{i_p}$ of $V^{\ot p}$
\be
A^{(p)}(R)_{i_1\dots
i_p}^{\;\;j_1\dots j_p} =
 u_{i_1\dots i_p}v^{j_1\dots j_p}.
\label{anti-s}
\ee
As one can show, the associated determinant $v$ possesses the property
$$
(x\ot v)\triangleleft R =  (v\ot x')\quad
{\rm  and}\quad (v\ot x)\triangleleft R = (x''\ot v)\quad \forall\, x\in
V
$$
where the correspondences $x\to x'$ and $x\to x''$ are some linear
maps $V\to V$. Demand them be scalar (multiple of the
identity map) and equal to each other. In other words, let
there exists a nonzero $a\in\K$ such that
$$
(x\ot v)\triangleleft R = a (v\ot x)\quad {\rm  and}\quad (v\ot
x)\triangleleft R =a (x\ot v).
$$
The Hecke symmetry $R$ satisfying such a requirement will be called
{\em admissible}.
In this case, by setting $\bar R =a^{-1} R$ we have
\be
(x\ot
v)\triangleleft{\bar R}= (v\ot x),\quad (v\ot
x)\triangleleft{\bar R} = (x\ot v). \label{r-bar}
\ee
The {\it  cancelling a column} is defined as a map $\psi$
$$
\lmp\stackrel{\psi}{\longrightarrow}\K:\quad \psi(v) = 1.
$$
Due to (\ref{r-bar}) $\psi$ obeys the following
condition
\be
\bar R\circ(\id\otimes\psi) = (\psi\otimes\id)\circ\bar R.
\label{rbar}
\ee
The above relation is valid also for the map $\psi^{-1}$ (inverse
to $\psi$). By definition $\psi$ and $\psi^{-1}$ are the morphisms
of the second kind. Any product (composition) $f\cdot g$ of morphisms of
the
both kinds gives {\em a categorical morphism}. Also, any tensor
product  $f\ot g$ of categorical morphisms will be a categorical
morphism by definition.

\begin{remark} {\rm
Condition (\ref{rbar}) (in a litlle bit more general form) is sometimes
included
in the system of axioms for braided categories. Thus, in \cite{T} a
braiding satisfying
such a condition is called {\it natural}. }
\end{remark}

Let us emphasize that for an object $\vl$ with a fixed embedding
$\vl\hookrightarrow V^{\ot m}$  a map $\vl\to \vl$ is a
categorical morphism if and only if it is a scalar map. This is the
reason
to call these objects {\em simple}. Such an observation plays an
important role in what follows. If a map $\xi:\,U\rightarrow W$
is proved to be a categorical morphism, then on each simple
component of $U$ it is a scalar map.

The category thus introduced is a monoidal and quasitensor one
 whose braidings are restrictions of
maps (\ref{braid}) onto the simple spaces $V_\lambda$ and their direct
sums.
 We call such a category {\em braided}. Note that its Grothendieck
(semi)ring is isomorphic to that of $sl(p)$-modules where
$p=\rk(R)$.

Moreover, if a Hecke symmetry $R$ is admissible then the corresponding
category $\B$ is rigid that is $\forall\,U\in {\rm Ob}(\B)$ the
dual spaces $U^*_\epsilon, \epsilon=r,l$ (right and left) are
also contained in ${\rm Ob}(\B)$. In particular, one can show that
$\lmpp$ is the dual (right and left) of $V$. This means that there
exist non-degenerate pairings
\be
\lmpp\ot V
\to \K\quad{\rm and}\quad V \ot \lmpp\to \K \label{l-r-pair}
\ee
which are categorical morphisms.

Let $\{x^i_r\}$ (resp. $\{x^i_l\}$) be the dual basis in
the right (resp. left) dual space to $V$ 
\be
<x_i,x^j_r>=\de_i^j,\quad <x_l^j,x_i>=\de_i^j.
\label{def:pair} 
\ee
The dual bases $\{x^i_r\}$ and $\{x^i_l\}$ can be expressed in the form 
\be
x^i_r=v^{a_1\dots\, a_{p-1}i}x_{a_1}\ot \dots\ot x_{a_{p-1}},\quad
x^i_l=v^{ia_1\dots\,a_{p-1}}x_{a_1}\ot \dots \ot x_{a_{p-1}},
\label{def:du-bas} 
\ee 
hence $\{x^i_{\epsilon}\}\in \Lambda_-^{p-1}(V),\, \epsilon=r,l$. Now, 
pairings (\ref{def:pair}) can be explicitly constructed by means of the 
categorical morphism $\psi$. As a consequence, pairings (\ref{def:pair}) 
become categorical morphisms justifying the identification of
$\Lambda_-^{p-1}(V)$ 
with the dual space of $V$ (see (\ref{l-r-pair})).

Introduce now a {\em categorical trace} which will play the central role
in all
our subsequent considerations. It is defined as a properly normalized
categorical
morphism $\tr_R:\,\End_\varepsilon(U) \to \K$. For details the reader is
referred
to \cite{GLS} and we only briefly outline this construction.

For any admissible Hecke symmetry R there exists an operator $Q$ that we
call
``inverse to $R$ by column'', i.e.,
$$
R_{ia}^{\;jb}\,Q_{bk}^{\;al}=\de_{i}^{\,l}\de_{k}^{\,j}\,\quad
\Leftrightarrow\quad
Q_{ia}^{\;jb}\,R_{bk}^{\;al}=\de_{i}^{\,l}\de_{k}^{\,j},
$$
where $R_{ia}^{\;jb}$ is the matrix of the Hecke symmetry in the
basis $x_i\otimes x_j$:
$$
(x_i\ot x_j)\triangleleft R = R_{ij}^{\;kl}\, x_k\ot x_l.
$$

Consider the matrices
$$
B_i^{\;j}=Q_{ai}^{\;aj},\quad C_i^{\;j}=Q_{ia}^{\;ja}.
$$
Evidently, they satisfy
\be
B_b^{\;a}R_{ai}^{\;bj}=\de_{i}^{\,j},\quad
R_{ia}^{\;jb}\,C_b^{\;a} =\de_{i}^{\,j}. \label{prop}
\ee
Extending these matrices to any objects $\vl$ in a proper way we
get the matrices $\bl$ and $\cl$ such that {\em the categorical
trace} ${\rm tr}_R$ on the space $\End_r(\vl)$ (resp.
$\End_l(\vl)$, see (\ref{endo})) is defined as follows
\be
{\mathrm{tr}}_R\,X=tr(\bl\cdot \hat X),\;\forall\,X\in\End_r(\vl)
\quad({\rm resp.}\;\;
{\mathrm{tr}}_R\, Y=tr( \cl\cdot\hat Y),\;\forall\,Y\in\End_l(\vl)).
\label{r-tr}
\ee
Here $tr$ is the usual matrix trace and $\hat X$ (resp. $\hat Y$ ) is
the matrix of the linear operator $V_\lambda\rightarrow V_\lambda$
corresponding to an element $X\in {\rm End}_r(V_\lambda)$ (resp.
$Y\in {\rm End}_l(V_\lambda)$).

The matrices $\bl$ and $\cl$ are
constructed in such a way that the map $X\rightarrow {\rm tr}_R X$
is a categorical morphism and, besides, the {\em  categorical
dimension}
$$
\dim_R(U)=
{\mathrm{tr}}_R\,\mathrm{id}_U\quad \forall\,U\in\Ob(\B)
$$
is an additive-multiplicative functional on the Grothendieck
(semi)ring. Then for any simple object $\vl$ one gets
$$
{\dim}_R(\vl)=\sll(q^{p-1},q^{p-3},...,q^{3-p},q^{1-p}),
$$
where $\sll$ is the Schur function in $p$ variables. A proof of this
fact is given in \cite{H}, \cite{GLS} (in another setting an
equivalent formula can be also found  in \cite{KW}).

A consequence of (\ref{def:du-bas}) is the fact that the maps \be
x^i_r\to x^a_lB_a^{\;i},\quad x^i_l\to x^a_r\,C_a^{\;i}
\label{l-r} \ee belong to $\Mor(\B)$. Therefore, the same is true
for the maps
$$
1\to  x^a_lB_a^{\;i}\, \ot x_i,\quad 1\to  x_i\ot x^a_r\,C_a^{\;i}
$$
since they are compositions of the map $\psi^{-1}$
$$
1 \stackrel{\psi^{-1}}{\longrightarrow} v = x_r^i\otimes x_i =
x_i\otimes x^{i}_l.
$$
and morphisms (\ref{l-r}).

Our next aim is to introduce some associative algebras
naturally connected to the categories involved. We consider
these algebras as braided analogs of the enveloping
algebras $U(gl(n))$ and $U(sl(n))$. Motivation will be
given later.

As a starting point of our construction we introduce elements
$l_i^{\;j}= x_i\ot x^j_r$ and form the matrix
\be
L=\|l_i^{\;j}\|, \quad 1\leq i,\,j\leq n=\ddim(V)
\label{def:l-mat}
\ee
where the lower index enumerates rows and
the upper one enumerates columns. Assume $R$ to be an admissible Hecke
symmetry and impose the following relations on the free algebra,
generated by all the elements $l_i^{\;j}$:
\be
R\,L_1\,R\,L_1-L_1\,R\,L_1\, R-\h(R\,L_1-L_1\,R)=0,\quad{\rm where}\quad
L_1=L\ot\id,\,\, \h\in\K \label{RE}
\ee
or, explicitly,
$$
R_{i_1i_2}^{\;a_1b_2}\,l_{a_1}^{\;b_1}\,R_{b_1b_2}^{\;c_1j_2}\,
l_{c_1}^{\;j_1} -
l_{i_1}^{\;a_1}\,R_{a_1i_2}^{\;b_1c_2}\,l_{b_1}^{\;c_1}\,
R_{c_1c_2}^{\;j_1j_2} - \h(R_{i_1i_2}^{\;j_1a}\,l_a^{\;j_2}-
l_{i_1}^{\;a}\,R_{a\,i_2}^{\;j_1j_2})=0.
$$
We call this relation the {\em  modified reflection
equation} (mRE) and the corresponding algebra, $\lhq$, {\em
the modified reflection equation algebra}. For
any Hecke symmetry $R$ with $q\not=1$ this algebra can be
obtained from the non-modified one (corresponding to
$\h=0$) by a shift of generators $l_{i}^{\;j}\rightarrow
l_{i}^{\;j} - a\,\de_i^{\,j}\,\id$ with
$a=\h(q-q^{-1})^{-1}$.

For the matrix $L$ we consider the maps
\be
L\to\id\quad
(l_i^{\;j}\to \de_i^{\,j})\quad {\rm and}\quad L\to L^{\ot
k}\quad (l_i^{\;j}\to l_i^{\;a_1}\ot l^{\;a_2}_{a_1} \ot \dots\ot
l^{\;j}_{a_{k-1}}) \label{morph}
\ee
which evidently belong to
$\Mor(\B)$. In what follows the  matrices $L^{\ot k}$ whose
entries are considered as elements of  the space $\lhq$ (or
$\slhq$ defined below) will be denoted $L^k$.

\begin{proposition}\label{prop:pi-1}
Let us set $l_i^{\;j}\triangleright x_k= x_iB_k^{\;j}$.
Then the image of the left hand side of (\ref{RE}) under
this map is equal to 0 if we put $\h=1$. Hence, we have a
representation
$$
\pi_1:\;\lhq\to \End_l(V),\;\h=1.
$$
\end{proposition}
{\bf Proof\ \ } Straightforward calculations. Suffice it to apply the
left hand side of (\ref{RE}) to an arbitrary element of $V$ and
use property (\ref{prop}) of matrix
$B$.\hfill\rule{6.5pt}{6.5pt}

\medskip

In $\End_l(V)$ we can choose the natural basis
$\{h_i^{\;j}\} =x_i\ot x^j_l$
$$
h_i^{\;j} \triangleright x_k=\de^{\,j}_k\, x_i
$$
and identify the element $l_i^{\;j}$ with $h_i^{\;a}B^{\;j}_a$.
Since matrix $B$ is non-degenerate, we can consider $\{l_i^{\;j}\}$ as
another basis in this space. In the basis $\{h_i^{\;j}\}$  the
representation
$\pi_1$  becomes tautological: $\pi_1(h_i^{\;j}) =h_i^{\;j}$.

It is worth emphasizing a difference between these bases. For the
set $\{h_i^{\;j}\}$, the product
$$
h_i^{\;j}\ot h_k^{\;l}\to \de_k^{\,j}\,h_i^{\;l}
$$
is a categorical morphism while for $\{l_i^{\;j}\}$  the maps
given in (\ref{morph}) are categorical morphisms.

Thus, we have identified the vector space Span$(l_i^j)$ with
$\End_l(V)$.
We denote this space by $gl_R(V)$ too and describe its subspace of
traceless
elements as follows.

For any rank $p$, the space $gl_R(V)$ decomposes into the
direct sum of two simple spaces, one of them is
one-dimensional. This one-dimensional component is
generated by the element $h_i^{\;i} = x_i\ot x^i_l =v$ or,
equivalently, by $\l = C_i^{\;j}\,l_j^{\;i}$. The elements
of other simple component of $\End_l(V)$ are called {\em
the traceless elements}. In the sequel such a space will be
denoted  $sl_R(V)$. Moreover, we define the algebra $\slhq$
as the quotient $\slhq=\lhq/\{\l\}$.

Propositions \ref{prop:pi-1} and \ref{prop:pi-m}, \ref{prop:slqh}
below suggest a new way of constructing the representation theory
of the algebras $\lhq$ and $\slhq$. In contrast with the usual
method valid in the case related to the quantum group when the
triangle decomposition of $L$ into the product of $L^+$ and $L^-$
is used, our approach works in the general setting (for arbitrary
admissible Hecke symmetry). In more detail this approach and the
proof of propositions \ref{prop:pi-m} and \ref{prop:slqh} will be
presented elsewhere.

Observe that all representations in question are
equivariant in the sense of the following definition.

\begin{definition}{\rm
Let A be either $\lhq$ or $\slhq$ and $U\in \B(V)$ be an object with
an associative product $U\ot U\to U$ which is a categorical morphism
(for
example, $U=\End_{\varepsilon}(W),\;\varepsilon=r,l$ or a direct
sum of tensor products of similar spaces).  We say that a map
$$
\pi_U:\; A\to U
$$
is {\it an equivariant  representation} if it is a representation (i.e.,
an
algebra morphism) and
its restriction to $gl_R(V)$ (resp. $sl_R(V)$) is a categorical
morphism.}
\end{definition}

Our next step is to define representations of $\lhq$ in all simple
spaces
$V_\lambda$  via a ``truncated coproduct" defined below. We constrain
ourselves
to the simplest case ${\rm rk}(R) = 2$. In this situation the simple
objects of
$\B$ are labelled by partitions of height 1: $\lambda = (m)$. The
corresponding
Young diagram has only one row of length $m$. For brevity, we will
write $V_{(m)}$
instead of $\vl,\,\, \la=(m)$.

Describe now the truncated coproduct. As we said above, for an
involutive
$R$, coproduct
(\ref{copr}) together with an antipode and a counit converts the
enveloping
algebra
into a braided Hopf algebra. Unfortunately, this coproduct is not ``Hopf
compatible''
with the algebraic structure of $\lhq$ at $q\not=1$. So, such a
coproduct
does not
allow us to  define maps $ \End_l(V)\to \End_l(V^{\ot m})$ which would
give
rise to
the higher representations of the algebra $\lhq$.

However, we define categorical morphisms which are
q-analogs of restricted maps
$$
\End_l(V)\to \End_l(V_{(m)})\hookrightarrow \End_l(V^{\ot m}).
$$
We will refer to the family of these categorical morphisms as
{\it the truncated coproduct}.

To construct them explicitly we take into account that
$V_{(m)}$ is the image of the projection $V^{\ot m}\to
\Lambda^m_{+}(V)$ where the corresponding projector $P_+^m$
is a polynomial in $R$. Then to an arbitrary element $X\in
\End_l(V)$ we assign the element $X_{(m)}\in
\End_l(V_{(m)})$ by the following rule \be
X_{(m)}\triangleright g=q^{1-m}\,[m]_q\, P_+^m
(X_{(1)}\triangleright g),\quad \forall\,g\in
\Lambda_+^m(V).\label{prol} \ee Here  $X_{(1)}=X\ot
\id_{(m-1)}$ and $[m]_q={q^m-q^{-m}\over q-q^{-1}}$ is the
q-analog of the integer $m$. Note that the map \be
\Delta_m:\;\End_l(V)\to \End_l(V_{(m)}),\quad X\to X_{(m)}
\label{map:pi-m} \ee is a categorical morphism due to the
structure of $P_+^m$. Composing $\pi_1$ with $\Delta_m$ we
get the map
$$
\pi_m:\; l_i^{\;j}\to \End_l(V_{(m)}).
$$

It is worth emphasizing that unlike $\Delta_m(l_i^{\;j})\in
\End_l(V_{(m)}) $ the elements $\pi_m(l_i^{\;j})$ are
considered to be operators (see section 3).

\begin{proposition}\label{prop:pi-m}
The image of the left hand side of (\ref{RE}) under the map $\pi_{m}$ is
equal to 0 at  $\h=1$.
So, we get a representation
$$
\pi_m:\; \lhq\to\End_l(V_{(m)}),\quad \h=1.
$$
\end{proposition}

\begin{remark}{\rm
If $\rk(R)=2$, it is easy to introduce a braided analog of
the Lie bracket in the space $sl_R(V)$. Taking into account
the decomposition
$$
sl_R(V)^{\ot 2}=V_{(4)}\oplus V_{(2)}\oplus V_{(0)}
$$
we set $[\;,\;]:\,V_{(4)}\oplus V_{(0)}\to 0$ and  require the map
$[\;,\;]:\,V_{(2)}\to sl_R(V)$ to
be a categorical morphism. This requirement defines the map
$[\;,\;]$ uniquely, up to a factor. This
bracket can be naturally extended to $gl_R(V)$ by the requirement
$[{\l},\, x]=0$ for any $x\in gl_R(V)$.
(A q-counterpart of the Lie algebra $sl(n)$ has been introduced in
\cite{LS}.)

Having introduced such a bracket, we can treat the algebra $\slhq$ as
the universal enveloping algebra
of the corresponding ``q-Lie algebra'' in a standard manner (the
parameter $\h$ depends on a normalization
of the bracket). Similarly, the algebra $\lhq$ can be treated as the
enveloping algebra of $gl_R(V)$.

However, we prefer to do without any ``q-Lie algebra''
structure. Similarly to the usual enveloping algebra, our
algebra $\slhq$ has the following properties. It is
generated by the space $sl_R(V)$ (more precisely, it is the
quotient of the algebra $T(sl_R(V))$ modulo an ideal
generated by some quadratic-linear terms). Moreover, its
representation theory resembles that of $U(sl(2))$  and,
being  constructed via the truncated coproduct, is
equivariant. This is the reason for considering the algebra
$\slhq$ as a proper ``braided analog'' of the enveloping
algebra $U(sl(2))$ (and similarly, the algebra $\lhq$ is
treated as the enveloping algebra of $gl_R(V)$).}
\end{remark}

The representations of the algebra $\slhq$ can be easily
deduced from those  of $\lhq$. To construct them, we set
$$
{\l}\triangleright x=0\,\,x\in V_{(m)}
$$
and preserve prolongation (\ref{prol}) for the elements of the
traceless component $sl_R(V)$.

\begin{proposition}\label{prop:slqh}
Thus defined maps are representations of the algebra $\slhq$ with
some $\h\not=0$.
\end{proposition}

We will refer to these representations of the algebra
$\slhq$ as {\it sl-representations} and  and keep the same
notation $\pi_m$ for them: $\pi_m:\,\slhq\to\End(V_{(m)})$.

The exact value of $\h$ in proposition \ref{prop:slqh} is
not important. Given a representation of $\slhq$  with some
$\h\not=0$, we can get a representation with another $\h$
renormalizing the generators in an appropriate  way. Note
that in  the $U_q(sl(2))$ case this method of constructing
representation theory of the algebra $\slhq$ was suggested
in \cite{DGR}.

Up to now we considered the ``left'' representations of the algebras
in question but we need also the ``right'' ones. Such
representations are given by appropriate maps
$$
l_i^{\;j}\to \End_r(V^*_{(m)}).
$$
Note that we do not specify which dual space --- right or left
--- is used in the formula above. In fact, it is of no importance
due to (iso)morphisms (\ref{l-r}). The right representation
$\bar\pi_m$ of $\lhq$ (and $\slhq$) in the space
$\Lambda^m_+(V^*_l)$ is introduced in the standard way as the map
$\bar\pi_m:\,\lhq\to \End_r(\Lambda^m_+ (V^*_l))$ given by the
formula
$$
<g\triangleleft\bar\pi_m(X),\, f>=<g,\,\pi_m(X)\triangleright f>
\quad {\rm for\  any}\quad f\in \Lambda^m_+(V),\; g\in
\Lambda^m_+(V^*_l).
$$

This construction is valid for an admissible Hecke symmetry of
any rank. The case $\rk(R)=2$ which we are dealing with leads to
additional technical simplifications. The point is that in this
case we can equip the space $V$ with a non-degenerate bilinear
form $V^{\ot 2} \rightarrow \Bbb K$ which is a categorical
morphism. This form allows us to identify $V$ with
$V^*_{\varepsilon},\;\varepsilon=r,l$ and, therefore, to define
representations $\bar\pi_m:\lhq\to\End_r(V_{(m)})$ and their
sl-counterparts $\bar\pi_m:\slhq\to\End_r(V_{(m)})$. These
representations are categorical morphisms also. Explicitly, such
a bilinear form can be written as follows \be <x_i,
\,x_j>=v^{-1}_{ij}.\label{pair} \ee
Here $\|v^{-1}_{ij}\|$ is the
matrix inverse to $\|v^{ij}\|$ which is invertible
as  has been shown in  \cite{G}.

Now, we define our main object --- the braided non-commutative sphere
--- as a quotient of $\slhq$.

\begin{definition} {\rm
Let  $\rk(R)=2$ and $\sigma\in\slhq$ be a nontrivial
quadratic central element (e.g., take ${\rm Tr}_R L^2$,
where ${\rm Tr}_R$ is defined in (\ref{Tr-R})). Fix
$\al\in\K$. The quotient $\slhq/\{\sigma- \al\}$ will be
called  {\it the braided non-commutative sphere}.}
\end{definition}

Here $\al$ is assumed to be generic. Below we consider
some polynomial (called Cayley-Hamilton) identities  whose
coefficients depend on $\al$. By demanding their roots to
be distinct we get more concrete restrictions on $\al$.

 In the particular case of the quantum sphere the
element $\sigma$ will be specified in section \ref{examp}.

Note that as in  the classical case this quotient has the
following spectral decomposition:
$$\slhq/\{\sigma-\al\}=\oplus_i V_{(2i)}.$$

\section{Index via braided Casimir element}

In \cite{GS1} we suggested a way to construct a family of
projective modules over the RE algebra (modified or not) by
means of the Cayley-Hamilton  identity. As was shown in
\cite{GPS}, the matrix $L$ satisfying (\ref{RE}) with any
even Hecke symmetry $R$ obeys to a polynomial relation \be
L^p+\sum_{i=0}^{p-1} \sigma_{p-i}(L)\, L^i=0, \quad
p=\rk(R), \label{CH} \ee where the coefficients
$\sigma_i(L)$ belong to the center $Z(\lhq)$ of the algebra
$\lhq$. This relation is called  {\it the Cayley-Hamilton
identity}.

Let us consider the quotient algebra $\lc=\lhq/\{I^{\chi}\}$,
where $\{I^{\chi}\}$ is the ideal generated by the elements
\be
z-\chi(z),\quad z\in Z(\lhq), \label{rea-xi}
\ee
where
$$
\chi:\; Z(\lhq)\to \K
$$
is a character of $Z(\lhq)$. After factorization to the algebra
$\lhq^{\chi}$ the coefficients in (\ref{CH}) become numerical
\be
L^p+\sum_{i=0}^{p-1} a_i\, L^i=0,\quad a_i=\chi(\sigma_{p-i}(L))
\label{chc-1}
\ee
(putting $\chi({\l})=0$ we obtain a quotient of $\slhq$ denoted as
$\slhq^\chi$).

Assuming the roots of the equation
$$
\mu^p+\sum_{i=0}^{p-1} a_i \mu^i=0
$$
to be distinct, one can introduce $p$ idempotents in the usual way
\be
e_i=\prod_{j\not=i} \frac{(L-\mu_j)}{(\mu_i-\mu_j)}, \quad
0\leq i\leq p-1. \label{proj}
\ee
If no character $\chi$ is fixed, then the roots $\mu_i$ can be treated
as elements of the
algebraic closure $\overline{Z(\lhq)}$ of the center $Z(\lhq)$
(or $\overline{Z(\slhq)}$
if $\chi({\l})=0$).

Besides the basic Cayley-Hamilton identity (\ref{CH}), we are interested
in the so-called {\em derived} ones which are valid for
some extensions of the matrix $L$ \cite{GS1}. A regular way
to introduce these extensions can be realized via  ``a
(split) braided Casimir element''. A particular case of
such a Casimir element corresponding to $q=1$ was used in
\cite{K}, \cite{R} in the study of the family algebras.

The {\it  braided Casimir element} is defined to be
\be 
\Cas=\sum_{i,j} l_i^{\;j}\ot h_j^{\;i} = \sum_{i,j}
l_i^{\;j}\ot l_j^{\;k}C_k^{\;i}. \label{def:cas} 
\ee 
Its crucial property is that the map
\be 
1\to \Cas \label{mor:cas} 
\ee 
belongs to $\Mor(\B)$.
Therefore, $\Cas$ is a central element of the
category $\B$ in the following sense
$$
(\Cas\ot U)\triangleleft {\overline R} = U \ot \Cas\quad
\forall\,U\in{\rm Ob}(\B)
$$
with ${\overline R}$ defined in the previous section. To prove this it
suffices
to observe that the above relation is obviously valid for the unity of
$\K$, hence,
for $\Cas$ due to the fact that (\ref{mor:cas}) is a categorical
morphism.

The braided Casimir element is a very useful tool in constructing
some extensions of the ``quantum matrix'' $L$. In order to
get the initial matrix $L$ as well as its higher analogs,
one should replace the elements $h_j^{\;i}\in \End_l(V)$ in
(\ref{def:cas}) with their images $\pi_m (h_j^{\;i})$
realizing them as matrices.
 For example, the matrix $L$ defined in (\ref{def:l-mat}) can be
represented
as follows
\be
L^t \stackrel{\mbox{\tiny def}}{=} \sum_{i,j} l_i^{\;j}\ot
\pi_1(h^{\;i}_j),
\label{ex:l}
\ee
where $L^t$ stands for the transposed matrix. The numerical matrix
$\pi_1(h^{\;i}_j)$
has the only nonzero $(i,j)$-th entry equal to 1.

Now, consider  the set of maps
$$
\Delta^{(2)}_m=\id\ot\Delta_m:\;\Cas \rightarrow
\lhq\ot \End_l(V_{(m)})\quad m=1,2\dots
$$
where $\Delta_m:\,\End_l(V)\to \End_l(V_{(m)})$ are defined by
(\ref{map:pi-m}), and another set
$$
\pi^{(2)}_m=\id\ot\pi_m:\;\Cas \rightarrow \lhq\ot
\Mat_{m+1}(\K)=\Mat_{m+1}(\lhq),\quad m=1,2\dots
$$
where the elements of $\End_l(V)$ are represented by the
corresponding matrices.

{\em An extension}
$L_{(m)}$ of the matrix $L$ is defined to be the image of $\Cas$ under
the map $\pi^{(2)}_m$
$$
L_{(m)}^t=\pi^{(2)}_m(\Cas) =\sum_{i,j} l_i^{\;j}\ot \pi_m
(h^{\;i}_j)\in \lhq\ot \Mat_{m+1}(\K)=\Mat_{m+1}(\lhq).
$$
An explicit example of such an extension related to the quantum
sphere will be given in section 4.

Below we will use the Cayley-Hamilton identity  for the
braided Casimir element which follows from (\ref{CH}); therefore,
we need  powers of $\Cas$. We define the product of
Casimir elements as the following composition of morphisms
$$
l_{i_1}^{\;j_1}\ot h_{j_1}^{\;i_1}\longrightarrow
l_{i_1}^{\;j_1}\ot 1\ot  h_{j_1}^{\;i_1}
\stackrel{(\ref{mor:cas})}{\longrightarrow}
l_{i_1}^{\;j_1}l_{i_2}^{\;j_2}\ot h_{j_2}^{\;i_2}h_{j_1}^{\;i_1}
\longrightarrow l_{i_1}^{\;a}l_{a}^{\;j_2}\ot h_{j_2}^{\;i_1}.
$$
That is \be \Cas^2=l_i^{\;a}l_a^{\;j}\ot h_j^{\;i},\quad
\Cas^3=l_i^{\;a}l_a^{\;b}l_b^{\;j}\ot h_j^{\;i},\quad {\rm etc}.
\label{pow-Cas} \ee

It is worth explaining the meaning of such a definition in
more detail. For matrices from ${\rm Mat}(\lhq)$ one
generally has $L^t\cdot L^t\not= (L^2)^t$. Therefore, a
care should be taken to define the product of the Casimir
elements in such a way as to preserve the proper
correspondence between the powers of $\Cas$ and those
of the matrix $L$. In particular, the definition might not
imply $\Cas^2\to (L^t)^2$. Our definition (\ref{pow-Cas})
is suitable in this sense because it gives the following
correspondence
$$
\Cas\to L^t, \quad \Cas^2\to (L^2)^t, \quad \Cas^3\to
(L^3)^t,\quad {\rm etc.}
$$
This means that if $L$ satisfies a Cayley-Hamilton
identity, then $\Cas$ obeys the same relation but
transposed as {\em a whole}, which  does not affect the
identity.

Definition (\ref{pow-Cas}) is directly transferred to the
product of extensions $L_{(m)}$ leading to the following
result \be (L_{(m)}^k)^t \stackrel{\mbox{\tiny def}}{=}
l_{i_1}^{\;j_1}\dots\,l_{i_k}^{\;j_k}\ot
\pi_m(h_{j_k}^{\;i_k})\dots\,\pi_m(h_{j_1}^{\;i_1}) \equiv
\left(\pi^{(2)}_m(\Cas)\right)^k, \label{Lm-pow} \ee where
the last equality is merely a conventional notation for the
preceding expression.

Define now a map ${\rm Tr}_R$ in the following way
\be
{\rm Tr}_R:\;\lhq\ot\Mat_{m+1}(\K)=\Mat_{m+1}(\lhq) \to \lhq,\quad
{\rm Tr}_R
\stackrel{\mbox{\tiny def}}{=} \id\ot{\rm tr}_R, \label{Tr-R}
\ee
where ${\rm tr}_R$ is the categorical trace (\ref{r-tr}) and the space
$\Mat_{m+1}(\K)$ is identified with $\End_l(V_{(m)})$.
 In particular, we have
$$
{\rm Tr}_R L ={\rm Tr}_R\, \pi^{(2)}_1 (\Cas)= l_i^{\;j}\ot {\rm
tr}_R (h_j^{\;i}) = l_i^{\;j}C_j^{\;i}=\l.
$$

In  the $U_q(sl(n))$ case the trace ${\rm Tr}_R L$
coincides with the {\em quantum trace} (cf. \cite{FRT})
which plays an important role in the theory of the RE
algebra.

Let us summarize the above construction once more. Given an
admissible Hecke symmetry $R$, we introduce the category $\B$ as
was shortly described in section \ref{sec:cat} and construct the
morphism (categorical trace) ${\rm tr}_R: \,\End_l(U)\to\K$,
$U\in {\rm Ob}(\B)$ which is defined by $R$. Then with the
category $\B$ we associate the algebra $\lhq$ defined by system
(\ref{RE}) and use this categorical trace in order to define the
map $\Tr_R$ sending
{\em matrices} with entries from $\lhq$  treated
as  elements of $\lhq\ot \End_l(V_{(m)})$ into the algebra
$\lhq$.

As we said above, the matrices  $L_{(m)}$ also satisfy
Cayley-Hamilton identities which we will call derived ones. 
Namely, there
exists a monic polynomial $\ch_{(m)}(t)$ of degree $m+1$
(recall that $\rk{R}=2$) whose coefficients belong to
$Z(\lhq)$ such that \be \ch_{(m)}(L_{(m)})=0,\quad
m=1,2,... \label{CH1} \ee

Pass now to the algebra $\lhq^{\chi}$ (see (\ref{rea-xi})) and
consider the image $\chc_{(m)}(t)$ of the polynomial
$\ch_{(m)}(t)$ in this algebra. Relation (\ref{CH1}) transforms
into a corresponding one in the algebra $\lhq^{\chi}$:
\be
\chc_{(m)}(L_{(m)})=0,\,\, m=1,2,... \label{CH2}
\ee
the coefficients of the polynomial $\chc_{(m)}(t)$ being numerical.
An explicit form of $\chc_{(m)}(L_{(m)})$ is determined by the
following proposition (the existence of identity
(\ref{CH1}) and the next proposition  will be proved in \cite{GS3}).

\begin{proposition} \label{conj}
Let $\mu_0$ and $\mu_1$ be roots of the polynomial
$\chc_{(1)}(t)$ (\ref{chc-1}) (at $p=2$ this polynomial is
quadratic). Then for each $m\ge 2$ the polynomial $\chc_{(m)}(t)$
is of the degree  $(m+1)$ and its roots $\mu_i(m)$
are given by the formula
\be
q^{m-1}\mu_{i}(m)=q^{-i}\,[m-i]_q\,\mu_0+ q^{i-m}[i]_q\,\mu_1+
[i]_q\,[m-i]_q\,\h,\quad i=0,1,...,m. \label{form}
\ee
\end{proposition}

Assuming the roots $\mu_i(m), \;0\leq i\leq m$, of the polynomial
$\chc_{(m)}(t)\;(m\ge 2)$ to be distinct we can introduce
idempotents $e_i(m)\in \lc\ot \End_l(V_{(m)})$ analogously to
(\ref{proj}) (to get uniform notations, we put $e_i=e_i(1)$).

If upon fixing some $m\ge 2$ one multiplies (\ref{CH2}) by
$L^n_{(m)},\; n\ge 0$, and then applies $\Tr_R$ to the resulting
equalities,
one obtains a recurrence  for  $\al_n(m)= \Tr_R L^n_{(m)},\; n\ge 0$.
The
general solution for such a recurrence is of the form
$$
\al_n(m)=\sum_{i=0}^m \mu_i^n(m)\, d_i(m),
$$
where $\mu_i(m)$ are the roots of the polynomial $\chc_{(m)}(t)$
(distinct by assumption) and the quantities $d_i(m)$  are defined
by the initial conditions, i.e., by the values $\Tr_R L^k_{(m)},\;
k=0,1,...,m$. Thus, we have  the following proposition.

\begin{proposition}\label{prop:9}
If the roots $\mu_i(m)$ of the polynomial $\chc_{(m)}(t)$ are all
distinct, then there exist $d_i(m)$ such that
$$
\Tr_R L^n_{(m)}=\sum_{i=0}^m \mu_i^n(m)\, d_i(m),\quad n=0,1,2,...
$$
\end{proposition}

The coefficients $d_i(m)$ in the above expansion are
functions in the roots $\mu_i(m)$. These functions  are
singular if there are coinciding roots. If we treat the
roots $\mu_i(m)$ as elements of $\overline{Z(\lhq)}$ (or
$\overline{Z(\slhq)}$), then the quantities $d_i(m)$
becomes elements of the field of fractions of this algebra.
 (Note that the case of multiple roots can be studied by
 the method of the paper \cite{DM}.)

Now, let us consider a representation $\opi_k:\,\lhq\to
\End_r(V_{(k)})$ of the algebra $\lhq$ (or $\slhq$) defined at
the end of section \ref{sec:cat}. It is easy to see that for a
generic $q$ the map $\opi_k$ is surjective and hence for any
$z\in Z(\lhq)$ the operator $\opi_k(z)$ is scalar. Therefore, we
can define a character $\chi_k:\, Z(\lhq)\to\K$ by taking
$$
\chi_k(z) = a_k(z),\quad{\rm where}
\quad \opi_k(z) = a_k(z)\,\id,\quad \forall\, z\in Z(\lhq)
$$
and denote $\ch_{k,m}(t)=\chc_{(m)}$  with $\chi=\chi_k$.
In what follows we will also use the notation \be
L_{(k,m)}^t=\opi_k(l_i^{\;j})\ot \pi_m(h_j^{\;i}).
\label{ltkm} \ee

We emphasize that  $L_{(k,m)}^t$ which is the image of the matrix
$L_{(m)}^t$ under the representation $\opi_k$ is treated as
an element of $\Mat_{m+1}(\Mat_{k+1}(\K))$. Also,
$L_{(k,m)}^t$ can be treated as an operator acting in the
space $V_{(k)}\ot V_{(m)}$. Indeed, if in the formula
(\ref{ltkm}) we consider $\opi_k(l_i^{\;j})$ and
$\pi_m(h_j^{\;i})$ as operators we get an operator acting
in the space $V_{(k)}\ot V_{(m)}$. More precisely, we put
the Casimir element $\Cas$ between the factors $V_{(k)}$ and
$V_{(m)}$ and apply it to these spaces via the
representations $\opi_k$ and $\pi_m$ respectively. This
operator generated by $\Cas$ and acting in the
product $V_{(k)}\ot V_{(m)}$ will be denoted $\Caskm$.

It is evident that the matrix $L_{(k,m)}$ satisfies the
Cayley-Hamilton identity
 \be
\ch_{k,m}(L_{(k,m)})=0 \label{raz} 
\ee 
which is a specialization of (\ref{CH2}) with $\chi=\chi_k$.
If the roots of the polynomial $\ch_{(k,m)}(t)$ are distinct, 
one can introduce idempotents $e_i(k,m)$ similarly to $e_i(m)$.

Applying the morphism
$$
\ttr=\tr_R^{(1)}\ot \tr_R^{(2)}:\;\; \End_r(V_{(k)})\ot
\End_l(V_{(m)})\to\K
$$
to all powers of the matrix $L_{(k,m)}$ and using the
Cayley-Hamilton identity for this matrix we can prove the
following proposition (similarly to proposition \ref{prop:9}).

\begin{proposition}
Let $\mu_i(k,m)$ be all the roots of the polynomial $\ch_{(k,m)}(t)$.
Let them
be distinct. Then there exist numbers $d_i(k,m),\;
0\leq i\leq m$  such that
$$
\ttr L_{(k,m)}^n= \sum_{i=0}^m \mu_i(k,m)^n\, d_i(k,m), \quad
n=0,1,2,...
$$
They are uniquely defined by the values of $\ttr L_{(k,m)}^l, \;
l=0,...,m$.
\end{proposition}

\begin{definition}{\rm
The quantities $\mu_i(m)$ and $d_i(m)$ (or $\mu_i(k,m)$ and
$d_i(k,m)$) will be called respectively {\it eigenvalues} and {\it
braided}
multiplicities of the matrix $L_{(m)}$ (or  $L_{(k,m)}$).}
\end{definition}

\begin{corollary} Let $f(t)$ be a polynomial (or a convergent series)
in $t$. Then
\begin{eqnarray*}
\Tr_R f(L_{(m)})\!\!\!&=&\!\!\!\sum f(\mu_i(m)) d_i(m),\\
\ttr f(L_{(k,m)})\!\!\!&=&\!\!\!\sum f(\mu_i(k,m)) d_i(k,m).
\end{eqnarray*}
\end{corollary}

In particular, taking as $f$ the polynomial in the right
hand side of (\ref{proj}) and its higher analogs we get the
following proposition.

\begin{proposition}
If the eigenvalues $\mu_i(m)$ (resp. $\mu_i(k,m)$) are distinct, then
\begin{eqnarray}
\Tr_R\, e_i(m)\!\!\!&=&\!\!\! d_i(m), \label{qind}\\
\ttr e_i(k,m)\!\!\!&=&\!\!\!d_i(k,m). \label{qInd}
\end{eqnarray}
\end{proposition}

\begin{definition}{\rm
The quantity $\ttr e_i(k,m)$ will be called {\it the q-index}
and denoted $\Ind(e_i(m),\opi_k)$.}
\end{definition}

\begin{remark}{\rm  Multiplying the trace by a factor results in a
modification
of the eigenvalues $\mu_i$ but does not affect the multiplicities
$d_i$. We are only interested in the latter quantities and
therefore can disregard a normalization of the trace. Similarly,
the multiplicities $d_i$ are stable under changes of the numeric
factor in (\ref{prol}). However, only for the factor $q^{1-m}[m]_q$ in
the definition of $X_{(m)}$ we get (\ref{form}).}
\end{remark}

\begin{remark}{\rm
On restricting to the $\uq$ case we emphasize that
similarly to \cite{K}, \cite{R} we deal with elements from
$(\lhq\ot\End(V_{(m)}))^{\uq}$, i.e., we consider
${\uq}$-invariant elements of this tensor product.
Representing $\lhq$ (or  $\slhq$) in the space $V_{(k)}$ we
obtain the space
$$
(\End(V_{(k)})\ot\End(V_{(m)}))^{\uq}.
$$
Thus, all our constructions are ``equivariant'' with
respect to the action of the quantum group (at least at the
$\uq$ case).}
\end{remark}

Taking into account that  $e_i(m)\in \lc\ot
\End_l(V_{(m)})$ and putting $\chi=\chi_k$ we get \be
e_i(k,m)=\opi^{(1)}_k(e_i(m)).\label{odin} \ee Finally, we
have \be \Ind(e_i(m),\opi_k)=\ttr e_i(k,m)=\ttr
\opi^{(1)}_k(e_i(m))=\tr_R\,\opi_k(\Tr_R\,
e_i(m)).\label{dva} \ee This justifies our treatment of the
quantity $\ttr e_i(k,m)$ as a braided (or q-)analog of the
Chern-Connes index. We would like to emphasize that
(\ref{odin}) and (\ref{dva}) are valid provided that the
eigenvalues $\mu_i(k,m)$ are pairwise distinct.

\begin{remark}{\rm
Speaking about the braided sphere, we are actually
dealing with a family of such spheres depending on the value of
the character $\chi=\chi_k$. So, if we treat the entries of the
idempotent $e_i(m)$ as elements of $\lc$ the q-index
$\Ind(e_i(m),\opi_k)$ is well-defined only for a special value of
$\chi$ depending on $k$.}
\end{remark}

Thus, due to (\ref{qInd}) and  (\ref{dva}) the computation
of our  q-index reduces to calculation of the braided
multiplicity $d_i(k,m)$. Now we will show how the latter
can be computed with the help of the operators $\Caskm$ defined 
above as images of the braided Casimir element:
$$
\Caskm:\;\;V_{(k)}\ot V_{(m)}\,\to\,V_{(k)}\ot V_{(m)}.
$$
Since map (\ref{mor:cas}) is a categorical morphism and since 
the representations $\opi_k$ and $\pi_m$  are equivariant we can 
conclude that each operator $\Caskm$ belongs to $\Mor(\B)$. This implies
that it is scalar on any simple component of the product
$V_{(k)}\ot V_{(m)}$.

Assuming $k\geq m$ one gets the following decomposition
$$
V_{(k)}\ot V_{(m)}=V_{(k+m)}\oplus V_{(k+m-2)}\oplus ... \oplus
V_{(k-m)}.
$$
Here we use that fact that the Grothendieck (semi)ring of the 
category in question is isomorphic to that of $sl(2)$-modules. 
Now, compute the trace of the operator $\Caskm$ using this 
decomposition. However, before doing it we would like to make 
the following observation.

Having fixed an object $U\in\Ob(\B)$, consider an arbitrary
linear operator ${\cal F}: U\to U$. What  is its trace? The
answer depends on the way this operator is realized. To
${\cal F}$, we can assign two elements: $F_l\in \End_l(U)$
and $F_r\in \End_r(U)$. In general, $\tr_R\,F_l\not=
\tr_R\,F_r$. However, under the additional assumption
${\cal F}\in\Mor(\B)$, the operator ${\cal F}$ should be
scalar on any simple component of $U$. For such an operator
the categorical trace is uniquely defined:
$$
\tr_R{\cal F} \stackrel{\mbox{\tiny def}}{=}\tr_R\,F_l =
\tr_R\,F_r.
$$
This follows from the trivial fact that $\tr_R\, \id$ is
the same for the right and left realization of the identity
operator. To sum up, if a linear operator belongs to ${\rm
Mor}(\B)$, then it is scalar on simple objects and its
categorical trace is uniquely defined. This observation
enables us to calculate $d_i(k,m)$.

Let $\mu_{i}$ be the eigenvalue of $\Caskm$ corresponding
to the component $V_{(k+m-2i)}$ ($0\leq i\leq m$). Then we
have \be \tr_R\,\Caskm^n=\sum_{i=0}^m \mu_i^n\,d_i,\quad
{\rm where}\quad d_i=\dim_R(V_{(k+m-2i)}).\label{op} \ee
Since $\dim_R(V_{(m)})=[m+1]_q$ in $\B$, we have the final
result.

\begin{proposition}\label{prop:19}
Let $k\geq m$ and let the eigenvalues $\mu_i(k,m)$ be all distinct.
Then we have
\be
\Ind(e_i(m),\opi_k)=[m+k-2i+1]_q,\,\, 0\leq i\leq m .\label{lab}
\ee
\end{proposition}
{\bf Proof\ \ } Under the hypothesis formulae (\ref{qind})--(\ref{dva})
are valid
and, therefore, the family of multiplicities $d_i(k,m)$ coincides with
that of $d_i$
from (\ref{op}). \hfill\rule{6.5pt}{6.5pt}
\medskip

In the next section we will see that if $k\geq m$ and $\opi_k$ are
sl-representations
in the $U_q(sl(2))$ case, then the eigenvalues $\mu_i(k,m)$ are
automatically distinct.

In the case of ``classical'' non-commutative sphere
(that is $q=1,\,\h\not=0$) we get
$$
\Ind(e_i(m),\opi_k)=m+k-2i+1
$$
which proves the formula given in
\cite{GS2}.

Observe that in Proposition \ref{prop:19} the eigenvalues $\mu_i(k,m)$
are numbered
according to decreasing (at $q=1$) dimensions of the components whereas
above their
numeration  was arbitrary.

\section{Example: quantum non-commutative sphere}
\label{examp}

Let us consider a particular case of the previous construction, namely,
the ``quantum non-com\-mu\-ta\-tive sphere''. In the framework of our
general 
approach we will introduce it using only the corresponding Hecke
symmetry, 
without any quantum group.

Let $V$  be a two dimensional vector space with a fixed base
$\{x_1,\, x_2\}$. Represent
the Hecke symmetry by the following matrix
$$
R=\left(\matrix{q&0&0&0\cr 0&\lambda&1&0\cr 0&1&0&0\cr
0&0&0&q}\right), \qquad \lambda = q-q^{-1}.
$$
The matrices $B$ and $C$ can be computed directly and after multiplying
by $q^2$
(which is just renormalization for the future convenience)
take the form
$$
B=\left(\matrix{q&0\cr 0&q^{-1}}\right),\quad
C=\left(\matrix{q^{-1}&0\cr 0&q}\right).
$$
We can choose the associated determinant as  $v = x_1\ot
x_2-q\,x_2\ot x_1$. Thus, we have
$$
\|v^{ij}\|=\left(\matrix{v^{11}&v^{12}\cr v^{21}&v^{22}}\right)=
\left(\matrix{0&1\cr -q&0}\right),\quad
\|v^{ij}\|^{-1}=\left(\matrix{0&-q^{-1}\cr 1&0}\right).
$$

Set
$$
l_1^{\;1} = a,\quad l_1^{\;2} = b,\quad l_2^{\;1} = c, \quad
l_2^{\;2} = d.
$$
In these generators the mRE algebra given by (\ref{RE}) takes the form
\be
\begin{array}{l@{\hspace{20mm}}l}
q ab - q^{-1}ba = \hbar b &q( bc - cb)  =(\lambda a -\hbar)(d-a)\\
q ca - q^{-1}ac  = \hbar c & q(cd - dc)  = c(\lambda a -\hbar)\\
ad-da=0 &  q(db - bd)   = (\lambda a -\hbar)b.
\end{array}
\label{ex:rea}
\ee
Represent the matrix $L$ in accordance with (\ref{ex:l})
$$
L^t=l_i^{\;j}\ot \pi_1(h^{\;i}_j) =
a\ot\pi_1(h_1^{\;1})+b\ot\pi_1(h_2^{\;1})+
c\ot\pi_1(h_1^{\;2})+d\ot\pi_1(h_2^{\;2})= \left(\matrix{a&c\cr
b&d}\right).
$$
(Recall that $\pi_1(h_i^{\;j})\triangleright x_k= \de^{\;j}_k
\,x_i$.) Taking (\ref{Tr-R}) into account we find
$$
{\l}=\Tr_R L =l_i^{\;j}C_j^{\;i} = q^{-1}\,a+q\,d.
$$

It is straightforward to check that ${\l}$ is a central element
in the mRE algebra. Now, let us consider the traceless
component $V_{(2)}=sl_R(V)$ of the space
$$
gl_R(V)=\span(a,b,c,d).
$$
For a basis in $sl_R(V)$ we take $\{b,\,c,\,g=a-d\}$.
Being reduced onto the traceless component of $gl_R(V)$,
system (\ref{ex:rea}) becomes
 \be
\begin{array}{c}
q^{2}gb - bg=\h (q+q^{-1}) b \\
gc - q^{2}cg= -\hbar (q+q^{-1}) c \\
(q^{2}+1)(bc-cb)+(q^{2}-1)g^{2}=\hbar (q+q^{-1})g.
\end{array}
\label{syst} \ee

Let us explicitly write the vector (two dimensional)
representations of $\slhq$  generated by (\ref{syst}).
Written  respectively in the bases $\{x_1,\,x_2\}$ and
$\{x^1_r,\,x^2_r\}$ the representations $\pi_1$ and
$\opi_1$ read on the generators:
$$
\pi_1(g)=\kappa\left(\matrix{q&0\cr 0&-q^{-1}}\right)\quad
\pi_1(b)=\kappa\left(\matrix{0&q^{-1}\cr 0&0}\right)\quad
\pi_1(c)=\kappa\left(\matrix{0&0\cr q &0}\right)\quad
\kappa\equiv \h\, \frac{q^2+1}{q^4+1}
$$
$$
\hspace*{-\coefkappa} \opi_1(g)=\kappa\left(\matrix{q&0\cr
0&-q^{-1}}\right)\quad \opi_1(b)=\kappa\left(\matrix{0&q\cr
0&0}\right)\quad \opi_1(c)=\kappa\left(\matrix{0&0\cr
q^{-1}&0}\right)\quad .
$$

In order to get the quantum non-commutative sphere we fix a
value of a nontrivial quadratic central element. As such an
element we take the coefficient $\sigma$ entering the
Cayley-Hamilton identity (\ref{CHn}). Then the {\em quantum
sphere} is obtained as the quotient of algebra (\ref{syst})
modulo the ideal $\{\sigma-\al\}$, for some $\alpha\in \K$.

An explicit form of the matrices $L$ and $L_{(2)}$ for
$\slhq$ is as follows. Taking the sl-representa\-tion
$\pi_2$ to construct $L_{(2)}$ we get
$$
L=L_{(1)}=\left(\matrix{
q[2]_q^{-1}g&b\cr
c&-q^{-1}[2]_q^{-1}g}
\right),
\qquad L_{(2)}=
q^{-1}\left(\matrix{
qg&[2]_qb&0\cr
q^{-1}c&(q-q^{-1})g&b\cr
0&q[2]_q c&-q^{-1}g}
\right)
$$
(the latter matrix is calculated in the basis $\{x_1^2,\,
qx_1x_2+x_2 x_1,\, x_2^2\}$).

One can directly check that the matrix $L$ satisfies the Cayley-Hamilton
identity of the form
\be
L^2-q^{-1}\,\h\,L + \sigma\, \id =0
\label{CHn}
\ee
where
\be
\sigma=-[2]_q^{-1}\,\Tr_R L^2 =
-[2]_q^{-1}\,([2]_q^{-1}g^2 + q^{-1}bc+qcb)\in
Z(\slhq).\label{cass}
\ee
The corresponding identity for the matrix $L_{(2)}$ reads
$$
L_{(2)}^3-2\h\, \frac{[2]_q}{q^2}\, L_{(2)}^2+
\frac{[2]_q^2}{q^2}\,(q^{-2}\h^2+\sigma)L_{(2)}
-\h\,\frac{[2]_q^3}{q^4}\,\sigma =0.
$$
This was shown in \cite{GS1} for a different normalization of $L_{(2)}$.

So, setting $\sigma=\al\in\K$ we come to the equation for $L$ with
numerical coefficients
$$
L^2-q^{-1}\,\h\, L+\al\,\id=0
$$
with the roots
$$
{\mu}_{0}={\mu}_{0}(1)=(q^{-1}\hbar-\sqrt{q^{-2}{\hbar}^{2}-4\alpha})/2,
\quad
{\mu}_{1}={\mu}_{1}(1)=(q^{-1}\hbar+\sqrt{q^{-2}{\hbar}^{2}-4\alpha})/2.
$$

The corresponding multiplicities (which coincide with $\Tr_R\, e_i(1)$
due to (\ref{qind})) are
$$
d_0(1)=\Tr_R\, e_0(1)=\Tr_R (L-\mu_1\,\id)(\mu_0-\mu_1)^{-1}=
\frac{[2]_q}{2}+\frac{[2]_q\hbar}{2\,\sqrt{{\hbar}^{2}-4\alpha
q^2}},
$$
$$
d_1(1)= \Tr_R\, e_1(1)=\Tr_R (L-\mu_0\,\id)(\mu_1-\mu_0)^{-1}=
\frac{[2]_q}{2}-\frac{[2]_q\hbar}{2\,\sqrt{{\hbar}^{2}-4\alpha
q^2}}.
$$

As for the matrix $L_{(2)}$, its eigenvalues can be found by means
of (\ref{form}) with $m=2$.

Our next aim is to compute the value of $\al$ corresponding to
the representation $\opi_k$ or, in other words, the value of
$\chi_k(\sigma)$. Clearly,  this value does not change
if we replace $\opi_k$ by $\pi_k$. Such a value (for a  Casimir
element being a multiple of (\ref{cass})) was computed in
\cite{DGR}. Using this result  we get
\be
\al=\chi_k(\sigma)=-\frac{\h^2\,[k]_q\,[k+2]_q}{q^2(\,[k+2]_q-[k]_q)^2}.
\label{sig-val}
\ee

This implies that
\be
\sqrt{q^{-2}{\hbar}^{2}-4\alpha}=\pm
\,{{q^{-1}\,[2]_q\, [k+1]_q}\over{[k+2]_q-[k]_q}}\h.\label{form1}
\ee
Choosing the positive sign in the right hand side of this formula we get
$$
\mu_0(k,1)={{-q^{-1}\h [k]_q}\over{[k+2]_q-[k]_q}},\qquad
\mu_1(k,1)={{q^{-1}\h [k+2]_q}\over{[k+2]_q-[k]_q}},
$$
$$
d_0(k,1)=[k+2]_q,\qquad d_1(k,1)=[k]_q.
$$
Note that the eigenvalues $\mu_i(k,1),\;i=0,1$, are distinct for
all $k\geq 1$.

\begin{proposition} On the quantum non-commutative sphere we have
$$
d_i(m)={{q^{m-2i}+q^{-m+2i}}\over{2}}+ {{[m-2i]_q
[2]_q\h}\over{2\sqrt{{\hbar}^{2}-4\alpha q^2}}},\,\,0\leq
i\leq m.
$$
\end{proposition}

{\bf Proof\ \ } Suffice it to check that
$$
\opi_k(d_i(m))\,[k+1]_q=[k+m+1-2i]_q.
$$
It can be easily done with the help of the following formula
$$
[k+m]_q+[k-m]_q=k_q\,(q^m+q^{-m}).
$$
\hfill\rule{6.5pt}{6.5pt}
\medskip

In fact, this proposition is valid for any braided sphere
since (\ref{sig-val}) can be shown to be true for any
admissible Hecke symmetry of rank 2. Also, note that for
$q=1$ we get formula (32) from \cite{K}.

\section{Concluding remarks}

1. If $R$ is an involutive symmetry, then there exists a natural
generalization of the notion of cyclic (co)homology, and categorical
trace becomes an ``$R$-cyclic cocycle''. However, we do not know any
generalization of this notion to the case of non-involutive $R$. Thus,
we do not know what should be a proper generalization of the ``highest
Chern-Connes indices'' to this case.

2. As we said in Introduction, in the classical case,
pairing (\ref{Ind}) reduces to (\ref{ind}) thanks to a
suitable definition of the group $K_0(A)$ based on the
notion of {\it stable isomorphic modules}. Two modules
$M_1$ and $M_2$ are called {\it stable isomorphic} (and are
identified in $K_0(A)$) if the corresponding idempotents
$e_1$ and $e_2$ (to be extended by 0 if necessary) become
``similar'': $e_1=Pe_2P^{-1}$ with some invertible  $P\in
\Mat_n(A)$  (cf. \cite{Ro} for detail).

However, in the braided case this equivalence implies neither
$\Tr_R\,e_1=\Tr_R\,e_2$ nor
$\ttr\,\pi(e_1)=\ttr\,\pi(e_2)$ for any representation $\pi$ of the
algebra
in question.
Let us point out that for the usual trace the latter relation is
satisfied because of the
particular property of the classical trace: it is stable under the
change
$M\to
PMP^{-1},\; M,P\in \Mat_n(\K)$, whereas the categorical trace is not.

Nevertheless, as follows from our computation (at least for the
quantum non-commutative sphere)
\be
\Tr_R\, e_i(m)=\Tr_R\,
e_{i+1}(m+2),\quad 0\leq i\leq m,\quad \forall\,m\ge 0.
\label{last}
\ee
In order to justify once more the use of the categorical
trace we want to emphasize that if we replace $\Tr_R$ by
the usual trace, then formula (\ref{last}) becomes wrong.

The modules for which the corresponding  idempotents have equal traces
$\Tr_R$  will be called {\em trace-equivalent}. It is easy to see that,
for
a generic $q$, the modules from the sequence
$$
e(0),\; e_0(m),\; e_{m}(m),\quad  m=1,2,...
$$
are not trace-equivalent. Thus, the set of classes of
trace-equivalent modules is labelled by $n=m-2i\in \Z$.
This looks like the Picard group of the usual sphere.

The problem whether the projective modules related to the
idempotents $e_{i}(m)$ and $e_{i+1}(m+2)$ are equivalent in
the conventional sense or not seems, however, to be
somewhat difficult even for $q=1$. Fortunately, the notion
of trace-equivalent modules suffices for our purposes. In
terms of $n=m-2i$ we can represent (\ref{lab}) as follows
\be 
\Ind(e_i(m),\opi_k)=[n+k+1]_q. \label{IInd} 
\ee 
Thus,
the q-index depends only on $n$ labelling classes of
trace-equivalent modules and $k$ labelling classes of
irreducible representations  of the algebra $\slhq$.

3. Consider a ``q-commutative'' analog of our algebras, i.e., set
$\h=0$.
In this case the algebra $\slhq$ does not have meaning of an enveloping
algebra
and we do not consider its representations. Now we take the classical
commutative
counterpart as a pattern.

Let us realize the usual sphere as a complex projective
variety. Then the line bundles ${\cal O}(n)$ and ${\cal
O}(-n),\,\,n\geq 0$, become analogs of our modules
corresponding to the idempotents $e_0(n),\,\,e_n(n)$. Which
line bundle corresponds to which projective module depends
on the holomorphic structure on the sphere (in our setting
the result depends on the sign of the root in
(\ref{form1})).

Let us consider the Euler characteristic
$$
\chi({\cal O}(n))=\dim\,H^0({\cal O}(n))-\dim\,H^1({\cal O}(n))
$$
of the bundle ${\cal O}(n)$ (for $n\geq 0$ it gives the
dimension of the space of global sections). Due to the
Riemann-Roch theorem we have $\chi({\cal O}(n))=n+1$ which
coincides with the above quantity  $[n+k+1]_q$ at $k=0$ and
$q=1$. So, for $q=1$ we consider index (\ref{IInd}) as a non-commutative
analog of the Euler characteristic of the class of the
idempotents  $e_i(m)$ with $n=m-2i$. Similarly, as the
q-analog of the Euler characteristic for a generic $q$ (in
the q-commutative case) we consider the specialization of
the $q$-Chern-Connes index at $k=0$, that is $[n+1]_q$.

4. As we said above, our quantum sphere is close to the Podles one but
it is defined without any quantum group as a quotient of the mRE
algebra. 
So, a question arises: whether it is possible to equip the RE algebra
(modified or not) or its quotient (the quantum sphere) with an
involution?
The answer is positive: it is easy
to see that the involution operator $*$ given by
$$
*b=c,\quad *c=b,\quad *g=g
$$
possesses the classical property $*(x\,y)=(*y)\,(*x)$ and is compatible
with (\ref{syst})
(here we assume $\h$ and $q$ to be real).

This fact is not surprising: such an involution exists for any
mRE algebra with the so-called real type $R$ (cf. \cite{M}). This
involution, considered as an operator
in $\End_\varepsilon(V)$, however, is not a categorical morphism since 
the
Euclidean pairing in the space $V$ is not.
The only (up to a factor) pairing in
$V$ which is a categorical morphism is given by (\ref{pair}). Using this
paring it is possible to define an involution in
$\End_\varepsilon(V)$ which does not satisfy the above classical
property (it looks like an involution in a super-algebra)  but we
do not need it at all (cf. \cite{DGR} for a discussion).

5. Let us mention the Poisson structures corresponding to
the quantum sphere and, more generally, to quantum orbits
related to the quantum group $\Uq$. On a generic orbit in
$sl(n)^*$, there exists a family of the so-called
Poisson-Lie structures (cf. \cite{DGS}). Their quantization
(in general, formal deformational) leads to algebras
covariant with respect to $\Uq$. But in this family only
one bracket (up to a numerical factor) is compatible with
the Kirillov one. Namely, the simultaneous quantization of
the corresponding ``Poisson pencil'' gives rise to the
quantum algebras which are appropriate quotients of $\slhq$
(``quantum orbits''). They depend on two parameters and the
particular case $\h=0$ is considered as ``q-commutative''
(the reader is referred to \cite{GS1} for  detail).

However, the properties of quantum algebras arising from
the Kirillov bracket alone and those arising from the above
pencil are different. The Kirillov  structure is symplectic
and for it there exists an invariant (Liouville) measure.
It gives rise to the classical trace in the corresponding
quantum algebra. On the contrary, the other brackets from
the Poisson pencil are not symplectic and they have no
invariant mesure. Their quantization leads to the algebras
with trace but this trace is braided. It is these algebras
and their ``non-quasiclassical'' analogs which are the main
objects of the present paper.

Also note that the  Poisson-Lie structures non-compatible with the
Kirillov bracket give rise to
one-parameter quantum algebras. It seems that for these algebras there
is
no reasonable means to
construct meaningful projective modules. On the usual sphere such
structures
do not exist
due to its low dimension.

6. The scheme presented in this paper is valid for the
``braided orbits'' related to the Hecke symmetries of
higher rank. Being  braided counterparts of semisimple but
not necessary generic orbits in $sl(n)^*$, such orbits can
be defined (at least for the $\Uq$ case) by methods of the
paper \cite{DM}. Thus, the ``easy part'', namely, the fact
that the q-index is nothing but a q-dimension of a
component in some tensor product, can be straightforwardly
generalized. The proof of an analog of formula (\ref{form})
is, however, much harder. Nevertheless, our low dimensional
computations make the following conjecture very plausible.

\begin{conjecture}
Let $\mu_i$ $1\le i\le p $ be roots of the polynomial
$\chc_{(1)}(t)$ (\ref{chc-1}). Then for $\forall\,m\ge 2$ the degree
of the polynomial $\chc_{(m)}(t)$ reads
$$
{\rm deg}(\chc_{(m)}(t)) =
{m+p-1\choose m}
$$
and its roots are given by the formula
\be
q^{m-1}\mu_{k_1\dots
k_p}(m)=\sum_{i=1}^p\frac{[k_i]_q}{q^{m-k_i}}\,\mu_i+
\xi_p(k_1,\dots,k_p)\,\hbar,\quad k_i\ge 0, \; \;k_1+\dots +k_p = m,
\ee
where $\xi_p(k_1,\dots,k_p)$ is the symmetric function defined as
follows
$$
\xi_p(k_1,\dots,k_p) =
\sum_{s=2}^{p}q^{k_1+k_2+\dots+k_s-m}[k_s]_q[k_1+k_2+\dots+k_{s-1}]_q.
$$
\end{conjecture}


\begin{thebibliography}{WWW}

\bibitem[C]{C}  A.Connes {\em Cyclic cohomology, Quantum group
symmetries
and Local Index formula for $\uq$}, QA/0209142.

\bibitem[CP]{CP} V.Chary, A.Pressley {\em A Guide to Quantum Groups},
Cambridge University Press, 1994.

\bibitem[DM]{DM} J.Donin, A.Mudrov {\em  Explicit equivariant
quantization
on (co)adjoint orbits of $GL(n, C)$}, QA/0206049.

\bibitem[DGR]{DGR} J.Donin, D. Gurevich, V.Rubtsov
{\em Quantum hyperboloid and braided modules},  Alg\`ebre
non-commutative, qroupes quantiques et invariants (Reims, 1995),
pp. 103--108,  S\'emin. Congr. 2, Soc. Math. France, Paris, 1997.

\bibitem[DGS]{DGS} J.Donin, D.Gurevich, S.Shnider
{\em Double quantization on some orbits in the coadjoint representations
of simple groups}, Comm. Math. Phys.
204 (1999), pp. 39--60.

\bibitem[FRT]{FRT} L.Faddeev, N.Reshetikhin, L.Takhtajan
{\em Quantization of Lie groups and Lie algebras}, Algebra
i Analiz, vol. 1, no 1 (1989), pp. 178--206 (in Russian);
English translation in: {\it Leningrad Math. J.} 1 (1990),
pp. 193--226.

\bibitem[G]{G} D.Gurevich {\em Algebraic aspects of the quantum
Yang-Baxter equation}, Leningrad Math. J. 2 (1991), pp. 801--828.

\bibitem[GLS]{GLS} D.Gurevich, R.Leclercq, P.Saponov
{\em Traces in braided categories}, J. Geom. Phys. 44
(2002), pp. 251--278.

\bibitem[GPS]{GPS} D.Gurevich, P.Pyatov, P.Saponov {\em Hecke
Symmetries and Characteristic Relations on Reflection Equation
Algebra}, Lett. Mat. Phys. 41 (1997) pp. 255--264.

\bibitem[GS1]{GS1} D.Gurevich, P.Saponov {\em Quantum line bundles via
Cayley-Hamilton identity}, J. Phys. A: Math. Gen.  34 (2001), pp.
4553 -- 4569.

\bibitem[GS2]{GS2} D.Gurevich, P.Saponov {\em Quantum line bundles on
noncommutative sphere}, J. Phys A: Math. Gen.  35 (2002),
pp. 9629 -- 9643.

\bibitem[GS3]{GS3} D.Gurevich, P.Saponov {\em
Geometry of noncommutative orbits }, in preparation.

\bibitem[Ha]{Ha} P.Hajac {\em Bundles over the quantum sphere
and noncommutative index theorem}, K-theory 21 (2001),
pp.141--150.

\bibitem[HM]{HM} P.Hajac, S.Majid {\em Projective module description of
the Q-monopole}, Comm. Math. Phys. 206 (1999), pp. 247--264.

\bibitem[H]{H} Phung Ho Hai {\em On matrix quantum groups of type
$A_n$}, Int. J. Math. 11 (2000), pp. 1115--1146.

\bibitem[KW]{KW} D.Kazhdan, H.Wenzl {\it Reconstructing Monoidal
Categories},
Adv. in Soviet. Math. 16, part 2 (1993), pp. 111--136.

\bibitem[K]{K} A.Kirillov {\em Introduction to family algebras},
Moscow Math.J. 1 (2001) pp.49--64.

\bibitem[L]{L} J-L.Loday {\em Cyclic homology}, Springer, 1998.

\bibitem[LS]{LS} V.Lyubashenko, A.Sudbery {\em Generalized Lie
algebras of $A_n$ type}, J. Math. Phys. 39 (1998), pp. 3487--3504.

\bibitem[M]{M} S.Majid {\em Foundations of quantum group theory},
Cambridge University Press, 1995.

\bibitem[MNW]{MNW} T.Masuda, Y.Najagami, J.Watanabe {\em Noncommutative
Differential Geometry on Quantum
Two Sphere of Podles. 1: An Algebraic Viewpoint},  K-theory 5 (1991),
pp. 151--175.

\bibitem[R]{R} N.Rozhkovskaya {\em Family algebras of representations
with simple spectrum}, preprint ESI-Vienna 1045 (2001).

\bibitem[Ro]{Ro} J.Rosenberg {\em Algebraic K-theory and its
application},
Springer-Verlag 1994.

\bibitem[P]{P} P.Podles {\em Quantum spheres}, Lett. Math. Phys. 14
(1987),
pp. 193-202.

\bibitem[T]{T} V.G.Turaev {\em Quantum Invariants of Knots and
3-Manifolds,} Walter de Gruyter, Berlin, NY 1994.

\end{thebibliography}
\end{document}